\documentclass{article}%
\usepackage{amsmath}
\usepackage{graphicx}%
\usepackage{amsfonts}%
\usepackage{amssymb}
\newtheorem{theorem}{Theorem}

\newtheorem{corollary}[theorem]{Corollary}

\newtheorem{proposition}[theorem]{Proposition}

\begin{document}
\bigskip

\bigskip

\bigskip

\bigskip

\bigskip

\bigskip

{\Huge Dynkin's isomorphism }

{\Huge without symmetry}

{\Huge \bigskip}

Yves Le Jan

Math\'{e}matiques

Universit\'{e} Paris 11.

91405 Orsay. France

\bigskip

\bigskip yves.lejan@math.u-psud.fr

\section{Introduction:}

The purpose of this note is to extend Dynkin isomorphim involving functionals
of the occupation field of a symmetric Markov processes with non polar points
and of the associated Gaussian field to a suitable class of non symmetric
Markov processes. This was briefly proposed in \cite{LJ2}\ using Grassmann
variables, extending to the non symmetric case some results of \cite{LJ1}.
Here we propose an alternative approach, not relying on Grassmann variables
that can be applied to the study of local times, in the spirit of \cite{MR}.
It works in general on a finite space and on an infinite space under some
assumption on the skew symmetric part of the generator which is checked on two examples.

\section{The finite case}

In this section, we will prove two formulas given in Proposition 1 and
Corollary 2 which relate the local time field of a non symmetric Markov
process on a finite space to the square of the associated complex ''twisted''
Gaussian field. The result and the proof appear to be a direct extension of
Dynkin's isomorphism.

\subsection{Dual processes}

Let us first consider the case of an irreducible Markov process on a finite
space $X$, with finite lifetime $\zeta$, generator $L$ and potential
$V=(-L)^{-1}$.

Let $m=\mu V$ for any nonnegative probability $\mu$ on $X$. $\operatorname{Re}%
$call that $L$\ can be written in the form $L=q(I-\Pi)$\ with $q$ positive and
$\Pi$ submarkovian . Then the $m$-adjoint $\widehat{L}$\ can be expressed
similarly with the same $q$\ and a possibly different submarkovian matrix
$\widehat{\Pi}$. Moreover, $m=\widehat{\mu}\widehat{V}$, with $\widehat{\mu}%
$\ the law of $x_{\zeta-}$ under $\mathbb{P}_{\mu}$.

Note that for any $z=x+iy\in\mathbb{C}^{X}$, the ''energy'' $<Lz,\overline
{z}>_{m}=\sum-(Lz)_{x}\overline{z}_{x}m_{x}.$ is nonnegative as it can be
written $\frac{1}{2}(\sum C_{x,y}(z_{x}-z_{y})(\overline{z_{x}}-\overline
{z_{y}})+<\Pi1+\widehat{\Pi}1-2,z\overline{z}>_{m})$, with $C_{x,y}%
=C_{y,x}=m_{x}q_{x}\Pi_{x,y}$.\ The highest eigenvector of $\frac{1}{2}%
(\Pi+\widehat{\Pi})$\ is nonnegative by the well known argument which shows
that the module contraction lowers the energy. and it follows from the strict
submarkovianity that the corresponding eigenvalue is strictly smaller than
$1$. Hence there is a ''mass gap'': For some positive $\varepsilon$, the
''energy'' $<-Lz,\overline{z}>_{m}$ dominates $\varepsilon<z,\overline{z}%
>_{m}$ for all $z$.

\subsection{A twisted Gaussian measure}

Then, although $L$ is not symmetric, an elementary computation (given in a
more general context in the following section) shows that for any $\chi
\in\mathbb{R}_{+}^{X}$, denoting $M_{\chi}$ the diagonal matrix with
coefficients given by $\chi$,
\[
\frac{1}{(2\pi)^{\left|  X\right|  }}\int(e^{-<z\overline{z},\chi>_{m}%
}e^{\frac{1}{2}<Lz,\overline{z}>_{m}}\Pi dx_{u}dy_{u}=\det(-M_{m}L+M_{\chi
m})^{-1}%
\]

. As a consequence, differentiating with respect to $\chi_{x}$,%
\[
\frac{1}{(2\pi)^{\left|  X\right|  }}\int z_{x}\overline{z_{x}}%
(e^{-<z\overline{z},\chi>_{m}}e^{\frac{1}{2}<Lz,\overline{z}>_{m}}\Pi
dx_{u}dy_{u}=\det(-M_{m}L+M_{\chi m})^{-1}\frac{1}{m_{x}}(-L+M_{\chi}%
)_{xx}^{-1}%
\]

In a similar way by perturbation using a non diagonal matrix, one obtains%
\[
\frac{1}{(2\pi)^{\left|  X\right|  }}\int z_{x}\overline{z_{y}}%
(e^{-<z\overline{z},\chi>_{m}}e^{\frac{1}{2}<Lz,\overline{z}>_{m}}\Pi
dx_{u}dy_{u}=\det(-M_{m}L+M_{\chi m})^{-1}\frac{1}{m_{y}}(-L+M_{\chi}%
)_{xy}^{-1}%
\]

But with the usual notations for Markov processes, setting $l_{t}^{x}=\int
_{0}^{t\wedge\zeta}1_{\{x_{s}=x\}}\frac{1}{m_{x_{s}}}ds$ and $l_{\zeta}%
^{x}=l^{x}$, we have%
\[
\frac{1}{m_{y}}(-L+M_{\chi})_{xy}^{-1}=\mathbb{E}_{x}(\int_{0}^{\zeta
}e^{-<\chi,l_{t}>_{m}}dl_{t}^{y})
\]

Defining the path measure $\mathbb{E}_{x,y}$\ by: $\mathbb{E}_{x}(\int
_{0}^{\zeta}G(x_{s},s\leq t)dl_{t}^{y})=\mathbb{E}_{x,y}(G)$ the above
relation writes%

\[
\frac{1}{m_{y}}(-L+M_{\chi})_{xy}^{-1}=\mathbb{E}_{x,y}(e^{-<\chi,l>_{m}})
\]

It follows that we have proved the following:

\begin{proposition}
For any continuous function $F$ on $\mathbb{R}_{+}^{X}$
\begin{align*}
(\ast)\;\int z_{x}\overline{z}_{y}F(z_{u}\overline{z}_{u},u &  \in
X)e^{\frac{1}{2}<Lz,\overline{z}>_{m}}\Pi dx_{u}dy_{u}\\
&  =\int\mathbb{E}_{x,y}(F(l^{u}+z_{u}\overline{z}_{u},u\in X))e^{\frac{1}%
{2}<Lz,\overline{z}>}\Pi dx_{u}dy_{u}%
\end{align*}
\end{proposition}

\subsection{Positivity}

It should be noted that setting $\rho_{u}=\frac{1}{2}z_{u}\overline{z}_{u}$
and $z_{x}=\sqrt{\rho_{x}/2}e^{i\theta_{x}}$the image on $\mathbb{R}_{+}^{X}%
$\ of the normalized complex measure $\nu_{X}=$ $\frac{1}{(2\pi)^{\left|
X\right|  }}\det(-M_{m}L)e^{\frac{1}{2}<Lz,\overline{z}>_{m}}\Pi dx_{u}dy_{u}$
\ by the map $z_{u}\rightarrow\rho_{u}$ is an infinitely divisible probability
distribution $Q$ on $\mathbb{R}_{+}^{X}$ with density $\frac{1}{(2\pi
)^{\left|  X\right|  }}\det(-M_{m}L)\int e^{<L\sqrt{\rho}e^{i\theta}%
,\sqrt{\rho}e^{-i\theta}>_{m}}\Pi d\theta_{u}$. Note that the positivity is
not a priori obvious when $L$ is not $m$-symmetric. This important fact
follows\ easily by considering the moment generating function $\Phi
(s)=\frac{\det(-L)}{\det(-L+M_{s})}=\det(I+(-L)^{-1}M_{s})^{-1}$ defined for
all $s$ with non negative coordinates, positive and analytic. The expansion in
power series\ around any $s$ (which appears for example in \cite{VJ}) is explicit:

$\frac{\Phi(s+h)}{\Phi(s)}=\det(I+(-L+M_{s})^{-1}M_{h})^{-1}=\exp(-\log
(\det(I+(-L+M_{s})^{-1}M_{h})))$

$=\exp(\sum\frac{(-1)^{k}}{k}Tr([(-L+M_{s})^{-1}M_{h}]^{k})$.

As $(-L+M_{s})^{-1}$ is nonnegative, it implies that $\Phi$ is completely
monotone as in this last expression, all coefficients of $h$-monomials of
order $n$ are of the sign of $(-1)^{n}$. \footnote{Completely monotone
functions in several variables were already used in \cite{Bo}.}

Note that the same argument works for fractional powers of $\Phi(s)$ which
shows the infinite divisibility. Let us incidentally mention it has been known
for a long time (cf \cite{VJ} ) that this expansion can be simplified further
in terms of permanents.

For $x=y$, the above proposition then yields the following:

\begin{corollary}
For any continuous function $F$ on $\mathbb{R}_{+}^{X}$%
\[
(\ast\ast)\;\int\rho_{x}F(\rho_{u},u\in X)Q(d\rho)=\int\mathbb{E}%
_{x,x}(F(l^{u}+\rho_{u},u\in X))Q(d\rho)
\]
\end{corollary}

Note that this last formula is also obtained in \cite{Eis} after a direct
definition of the measure $Q$.

Another interpretation of this positivity and of infinite divisibility can be
given in terms of a Poisson process of loops. It will be developed in a
forthcoming paper but let us simply mention that $Q$ appears to be equal to
the distribution of the occupation field associated with the Poisson process
of loops canonically defined by the Markov chain.

\bigskip

REMARK: If $Y\subset X$, it is well known that the trace of the process on $Y$
is a Markov process the potential of which is the restriction of $V$ to
$Y\times Y$. The distribution $\nu_{X}$\ induces $\nu_{Y}$. Therefore the
formulas (*)\ and (**) on $X$ and $Y$ are consistent.

EXAMPLE: let us consider, as an example, the case where $X=\{1,2....N\}$,
$q_{i}=1$, $\Pi_{i,j}=\mathbf{1}_{i<N,j=i+1}$, $\mu_{i}=\mathbf{1}_{i=1}$,
$m_{i}=1$ and $(-L)_{i,j}^{-1}=\mathbf{1}_{i\leq j}$..

The characteristic polynomial of $\frac{1}{2}(\Pi+\widehat{\Pi})$\ is
$(-\lambda+\sqrt{\lambda^{2}-1})^{N}+(-\lambda-\sqrt{\lambda^{2}-1})^{N}$
hence one gets easily that the mass gap equals $2\sin^{2}(\frac{\pi}{2N})$

Under $\mathbb{E}_{x,x}$\ all local times vanish except $l^{x}$ which follows
an exponential distribution. Moreover an easy calculation shows that $Q$
reduces to a product of exponential distributions. The formula (**) reduces to
the convolution of two exponentials.

\section{The infinite case}

We now explain how in certain situations, the above can be extended to a
Markov process on an infinite space $X$. Of course, a Markov process for which
points are not polar can always be viewed elementarily as a consistent system
of processes on finite subspaces but we aim at a stronger representation
allowing to consider any functional of the occupation field. There are some
obvious obstructions to a generalization. The mass gap property does not
always hold: consider for example the case of a constant drift on an interval,
analogous to the above example. Some conditions have to be assumed in order
that the energy controls the skew-symmetric part of the generator.

\subsection{Some calculations in Gaussian space}

Let $H$ be a real Hilbert space with scalar product $<,>$. At first the reader
may suppose it finite dimensional and then check that the assumptions we will
make allow to extend the results to the infinite dimensional case.

Let $\phi$ be the canonical Gaussian field indexed by $H$. Given any ONB
$e_{k}$ of $H$, $w_{k}=\phi(e_{k})$ are independent normal variables. Recall
that for all $f\in H$, $\phi(f)=\sum_{k}<f,e_{k}>w_{k}$ and $E(e^{i\phi
(f)})=e^{-\frac{\left\|  f\right\|  ^{2}}{2}}$.

In the following, $\phi(f)$ can be denoted by $<\phi,f>$ though of course
$\phi$ does not belong to $H$ in general.

Let $K$ be any Hilbert-Schmidt operator on $H$. Note that $K\phi=\sum_{k}%
w_{k}Ke_{k}$ is well defined as a $H$-valued random variable, and that
$E(\left\|  K\phi\right\|  ^{2})=Tr(KK^{\ast})$.

Let $C$ be a symmetric non negative trace-class linear operator on $H$. Recall
that the positive integrable random variable $<C\phi,\phi>\in\mathbb{L}^{1}$
can be defined by $\sum_{k}<Ce_{k},e_{k}>w_{k}^{2}$ for any ONB diagonalizing
$C$.

Moreover, $E(e^{-\frac{1}{2}<C\phi,\phi>+i\phi(f)})=\det(I+C)^{-\frac{1}{2}%
}e^{-\frac{<(I+C)^{-1}f,f>}{2}}$ (the determinant can be defined as
$\prod(1+\lambda_{i})$, where the $\lambda_{i}$ are the eigenvalues of $C$.

In fact $\det(I+T)$ is well defined for any trace class operator $T$ as
$1+\sum_{n=1}^{\infty}Tr(T^{\wedge n})$ (Cf \cite{Sim}\ Chapter 3). It extends
continuously the determinant defined with finite ranks operators and it
verifies the identity:

$\det(I+T_{1}+T_{2}+T_{1}T_{2})=\det(I+T_{1})\det(I+T_{2})$. By Lidskii's
\ theorem, it is also given by the product $\prod(1+\lambda_{i})$ \ defined by
the eigenvalues of the trace class (hence compact) operator $T$, counted with
their algebraic multiplicity.

Let $\phi_{1}$ and $\phi_{2}$ be two independent copies of the canonical
Gaussian process indexed by $H$. Let $B$ be a skew-symmetric Hilbert-Schmidt
operator on $H$. Note that $<B\phi_{1},\phi_{2}>=\sum<Be_{k},e_{l}>w_{k}%
^{1}w_{l}^{2}=-<B\phi_{2},\phi_{1}>$ is well defined in $\mathbb{L}^{2}$\ and
that $E(e^{i<B\phi_{1},\phi_{2}>})=E(e^{-\frac{1}{2}\left\|  B\phi
_{1}\right\|  ^{2}})=\det(I+BB^{\ast})^{-\frac{1}{2}}$.

As $B$ is Hilbert-Schmidt, $BB^{\ast}$ is trace-class. The renormalized
determinant $\det_{2}(I+B)=\det((I+B)e^{-B})$ is well defined (Cf \cite{Sim}),
and as the eigenvalues of $B$ are purely imaginary and pairwise conjugated, it
is strictly positive. Moreover since $B^{\ast}=-B$, $\det_{2}(I+B)=\det
_{2}(I-B)$= $\det(I+BB^{\ast})^{\frac{1}{2}}$

Finally, it comes that $E(e^{i<B\phi_{1},\phi_{2}>})=\det_{2}(I+B)^{-1}$

More generally, setting $\psi=\phi_{1}+i\phi_{2}$,

$E(e^{-\frac{1}{2}<(C-B)\psi,\overline{\psi}>})=\det_{2}(I+C+B)^{-1}\exp(-Tr(C))$

(Recall that when $T$ is trace class, $\det_{2}(I+T)\exp(Tr(T))=\det(I+T)$).

Indeed, $E(e^{-\frac{1}{2}<(C-B)\psi,\overline{\psi}>})=E(e^{-\frac{1}%
{2}<C\phi_{1},\phi_{1}>-\frac{1}{2}<C\phi_{2},\phi_{2}>+i<B\phi_{1},\phi_{2}>})$

$=\det(I+C)^{-\frac{1}{2}}E(e^{-\frac{1}{2}(<(I+C)^{-1}B\phi_{1},B\phi
_{1}>+<C\phi_{1},\phi_{1}>)})$ (by integration in $\phi_{2}$)

$=\det(I+C)^{-\frac{1}{2}}\det(I+C-B(I+C)^{-1}B)^{-\frac{1}{2}}$

$=\det(I-(I+C)^{-\frac{1}{2}}B(I+C)^{-1}B(I+C)^{-\frac{1}{2}})^{-\frac{1}{2}%
}\det(I+C)^{-1}$

=$\det((I+(I+C)^{-\frac{1}{2}}B(I+C)^{-\frac{1}{2}})(I-(I+C)^{-\frac{1}{2}%
}B(I+C)^{-\frac{1}{2}}))^{-\frac{1}{2}}\det(I+C)^{-1}$

$=\det_{2}(I+(I+C)^{-\frac{1}{2}}B(I+C)^{-\frac{1}{2}})^{-1}\det(I+C)^{-1}%
$(since $\det_{2}(I+(I+C)^{-\frac{1}{2}}B(I+C)^{-\frac{1}{2}})=\det
_{2}(I-(I+C)^{-\frac{1}{2}}B(I+C)^{-\frac{1}{2}})\ $\ by skew symmetry as before)

$=\det_{2}(I+C+B)^{-1}\exp(-Tr(C))$.

\bigskip

Note that $I+C+B$ \ is always invertible, as $C+B$ is a compact operator and
$-1$ is not an eigenvalue.

\bigskip

Let $f_{1}$ and $f_{2}$ be two elements of $H$. Set $D(f)=<f,f_{1}>f_{2}$. For
small enough $\varepsilon$, $E(e^{-\frac{1}{2}(<(C-B)\psi,\overline{\psi
}>+\varepsilon\psi(f_{1})\overline{\psi}(f_{2}))})$

$=\det_{2}(I+C+B+\varepsilon D)^{-1}\exp(-Tr(C+\varepsilon D))$

$=\det_{2}(I+C+B)^{-1}\exp(-Tr(C))\det(I+\varepsilon(I+C+B)^{-1}D)^{-1}$.

Hence, differentiating both members at $\varepsilon=0$,
\[
E(\psi(f_{1})\overline{\psi}(f_{2})e^{-\frac{1}{2}(<(C-B)\psi,\overline{\psi
}>)})={\det}_{2}(I+C+B)^{-1}\exp(-Tr(C))Tr((I+C+B)^{-1}D)
\]
Therefore%
\[
\frac{E(\psi(f_{1})\overline{\psi}(f_{2})e^{-\frac{1}{2}(<(C-B)\psi
,\overline{\psi}>)})}{E(e^{-\frac{1}{2}(<(C-B)\psi,\overline{\psi}>)}%
}=<(I+C+B)^{-1}(f_{1}),f_{2}>.
\]
If $C$ is only Hilbert-Schmidt, we can consider only the renormalized ''Wick
square'' $:<C\phi,\phi>:=\sum_{k}<Ce_{k},e_{k}>(w_{k}^{2}-1)$ for any ONB
diagonalizing $C$ and $E(e^{-\frac{1}{2}:<C\phi,\phi>:+i\phi(f)})=\det
_{2}(I+C)^{-\frac{1}{2}}e^{-\frac{<(I+C)^{-1}f,f>}{2}}$

The results given above extend immediately as follows:%
\[
E(e^{-\frac{1}{2}:<(C-B)\psi,\overline{\psi}>:})={\det}_{2}(I+C+B)^{-1}%
\]%
\[
\frac{E(\psi(f_{1})\overline{\psi}(f_{2})e^{-\frac{1}{2}(:<(C-B)\psi
,\overline{\psi}>:)})}{E(e^{-\frac{1}{2}(:<(C-B)\psi,\overline{\psi}>:)}%
}=<(I+C+B)^{-1}(f_{1}),f_{2}>.
\]

\subsection{A class of Markov processes in duality}

Let $(V_{\alpha},\alpha\geq0)$ and $\widehat{V}_{\alpha}$ be two Markovian or
submarkovian resolvents in duality in a space $\mathbb{L}^{2}(X,\mathcal{B}%
,m)$ with generators $L$ and $\widehat{L}$ , such that:

1) Denoting $\mathcal{D=D}(L)\cap\mathcal{D}(\widehat{L})$, $L(\mathcal{D})$
\ is dense in $\mathbb{L}^{2}(m)$

2) $<-Lf,f>_{m}\geq$ $\varepsilon<f,f>_{m}$ for some $\varepsilon>0$ and any
$f\in\mathcal{D}$ (i.e. we assume the existence of a spectral gap: it can
always be obtained by adding a negative constant to $L$.

Let $H$ be the completion of $\mathcal{D}$ with respect to the energy norm. It
is a functional space space imbedded in $\mathbb{L}^{2}(X,\mathcal{B},m)$. Let
$A$ be the associated self adjoint generator so that $H=\mathcal{D}(\sqrt
{-A})$. On $\mathcal{D}$, $\frac{1}{2}(L+\widehat{L})=A$.

The final assumption is crucial to allow the control of the antisymmetric part:

3) $B=(-A)^{-1}\frac{L-\widehat{L}}{2}$ is a Hilbert-Schmidt operator on $H$.

Equivalently, $(-A)^{-\frac{1}{2}}\frac{L-\widehat{L}}{2}(-A)^{-\frac{1}{2}}$
is an antisymmetric Hilbert Schmidt operator on \ $\mathbb{L}^{2}(m)$ since
for any ONB $e_{k}$\ of $H$, $(-A)^{\frac{1}{2}}e_{k}$ is an ONB of
$\mathbb{L}^{2}(m)$. Note that $I-B$ is bounded and invertible on $H$ and that
$V_{0}=(I-B)^{-1}(-A)^{-1}$ maps $\mathbb{L}^{2}(m)$ into $H$. Indeed, one can
see first that on $\mathcal{D}$, $A(I-B)=L$ so that on $L(\mathcal{D)}$,
$V_{0}=(I-B)^{-1}(-A)^{-1}$

\bigskip

EXAMPLES

This applies to the case of the finite space considered above.

\ Let us mention other examples:

1) Diffusion with drift on the circle: $X=S^{1}$, $A=\frac{\partial^{2}%
}{\partial\theta^{2}}-\varepsilon$, $L-\widehat{L}=b(\theta)\frac{\partial
}{\partial\theta}$, where $b$ is a bounded function on $S^{1}$.

Indeed, considering the orthonormal basis $e^{ik\theta}$ in $\mathbb{L}%
^{2}(d\theta)$, $\frac{1}{\sqrt{k^{2}+\varepsilon}}e^{ik\theta}$ is an
orthonormal basis in $H=H^{1}$, and

$\sum_{k}\left\|  (-A)^{-\frac{1}{2}}b(\theta)\frac{\partial}{\partial\theta
}\frac{1}{\sqrt{k^{2}+\varepsilon}}e^{ik\theta}\right\|  _{\mathbb{L}%
^{2}(d\theta)}^{2}=\sum_{k,l}\frac{k^{2}}{k^{2}+\varepsilon}(\widehat
{b}(l-k))^{2}\frac{1}{l^{2}+\varepsilon}<\infty$

2) Levy processes on the circle: The Fourier coefficients $a_{k}+ib_{k}$ of
$L$ should verify $\sum_{k}(\frac{b_{k}}{a_{k}})^{2}<\infty$

\subsection{An extension of Dynkin's isomorphism}

Assume $X$ is locally compact and separable, and that functions of $H$ are
continuous. By the Banach-Steinhaus theorem, given any point $x\in X$, there
exists an element of $H$, denoted $\eta_{x}$ defined by the identity:
$f(x)=<\eta_{x},f>_{H}.$ Note that $\eta_{x}=\sum e_{k}(x)e_{k}$ for any
orthonormal basis of $H$.

The resolvent $V_{\lambda}$\ is necessarily Fellerian and induces a strong
Markov process. Denote by $l_{t}^{x}$ the local time at $x$ of this Markov
process. Let $x$ and $y$ be two points of $X$. Set: $<\eta_{x},\eta_{y}%
>_{H}=\sum e_{k}(x)e_{k}(y)=K(x,y)$ so that $A^{-1}f(x)=\int K(x,y)f(y)m(dy)$
or $A^{-1}f=\int\eta_{y}f(y)m(dy).$

Set $V_{0}(x,y)=<I+B)^{-1}\eta_{x},\eta_{y}>_{H}$ and note that $V_{0}(x,y)$
is a kernel for $V_{0}$. Indeed, for any $f,g\in\mathbb{L}^{2}(m)$,
$<V_{0}f,g>_{\mathbb{L}^{2}(m)}=<(I-B)^{-1}(-A)^{-1}f,g>_{\mathbb{L}^{2}%
(m)}=<(I-B)^{-1}(-A)^{-1}f,A^{-1}g>_{H}=\int f(x)g(y)V_{0}(x,y)m(dx)m(dy)$

As a consequence, $V_{0}(x,y)=\mathbb{E}_{x}(l_{\zeta}^{y})$.

Applying the construction of the section 3-1, we see the kernel $K(x,y)$ is
the covariance of a Gaussian process $(Z_{x}=\psi(\eta_{x})=\sum e_{k}%
(x)w_{k},\,x\in X)$, for any ONB $e_{k}$ of $H$.

More generally, for any non-negative finitely supported measure $\chi=\sum
_{1}^{N}p_{j}\delta_{u_{j}}$ on $X$, letting $C$ be the finite rank
operator$:C=\sum_{1}^{N}p_{j}\eta_{u_{j}}\otimes\eta_{u_{j}},$ set $V_{\chi
}(x,y)=<(I-B+C)^{-1}\eta_{x},\eta_{y}>_{H}$

In a similar way as above for $V_{0}$, we have $V_{\chi}(x,y)=\mathbb{E}%
_{x}(e^{-\int l_{t}^{z}\,\chi(dz)}dl_{t}^{y})$. Then, from section 3-1

$E(Z_{x}\overline{Z}_{y})e^{-\frac{1}{2}(<(C-B)\psi,\overline{\psi}>_{H}%
)})=E(e^{-\frac{1}{2}(<(C-B)\psi,\overline{\psi}>_{H})})<(I-B+C)^{-1}\eta
_{x},\eta_{y}>_{H}$

$=E(e^{-\frac{1}{2}(<(C-B)\psi,\overline{\psi}>_{H})})V_{\chi}(x,y)$

But $(<C\psi,\overline{\psi}>_{H}=\int\psi(\eta_{u})\overline{\psi}(\eta
_{u})\chi(du)=\int Z_{u}\overline{Z}_{u}\chi(du)$ and $<B\psi,\overline{\psi
}>_{H}=<(-A)^{-1}\frac{L-\widehat{L}}{2}\psi,\overline{\psi}>_{H}%
=\sum<(-A)^{-1}\frac{L-\widehat{L}}{2}e_{k},e_{l}>_{H}w_{k}^{1}w_{l}^{2}$.

On the other hand, at least formally in general but exactly in the finite
dimensional case,

$\int\frac{L-\widehat{L}}{2}Z_{u}\overline{Z}_{u}m(du)=<(-A)^{-1}%
\frac{L-\widehat{L}}{2}Z,\overline{Z}>_{H}=\sum w_{k}^{1}w_{l}^{2}%
<(-A)^{-1}\frac{L-\widehat{L}}{2}e_{k},e_{l}>_{H}$

Hence we can denote: $<B\psi,\overline{\psi}>_{H}$ by $<\frac{L-\widehat{L}%
}{2}Z,\overline{Z}>_{\mathbb{L}^{2}(m)}$. Therefore

$E(Z_{x}\overline{Z}_{y}e^{\frac{1}{2}<\frac{L-\widehat{L}}{2}Z,\overline
{Z}>_{\mathbb{L}^{2}(m)}-\int Z_{u}\overline{Z}_{u}\chi(du)})$

$=E\otimes\mathbb{E}_{x}(\int_{0}^{\zeta}e^{\frac{1}{2}<\frac{L-\widehat{L}%
}{2}Z,\overline{Z}>_{\mathbb{L}^{2}(m)}}e^{-\int l_{t}^{z}\,\chi(dz)-\int
Z_{u}\overline{Z}_{u}\chi(du)}dl_{t}^{y})$

Let $\mathbb{E}_{x,y}^{t}$\ denote the non-normalized law of the bridge of
duration $t$ from $x$\ to $y$. Set $\mathbb{\mu}_{x,y}=\int_{0}^{\infty
}\mathbb{E}_{x,y}^{t}dt$, so that $\mathbb{\mu}_{x,y}(1)=$ $\mathbb{E}%
_{x}(l_{\zeta}^{y})=V_{0}(x,y)$ and denote $l^{z}$ the local time at $z$ for
the duration of the bridge.

Then $E(Z_{x}\overline{Z}_{y}e^{\frac{1}{2}<\frac{L-\widehat{L}}{2}%
Z,\overline{Z}>_{\mathbb{L}^{2}(m)}-\int Z_{u}\overline{Z}_{u}\chi(du)})$

$=E(e^{\frac{1}{2}<\frac{L-\widehat{L}}{2}Z,\overline{Z}>_{\mathbb{L}^{2}%
(m)}-\int Z_{u}\overline{Z}_{u}\chi(du)}\int e^{-\int l^{z}\,\chi(dz)}%
d\mu_{x,y})$

Finally, we get that for any continuous bounded function $F$ of $N$ non
negative real coordinates, and any $N$-uple of points $u_{j}$ in $X$,%

\[
E(Z_{x}\overline{Z}_{y}e^{\frac{1}{2}<\frac{L-\widehat{L}}{2}Z,\overline
{Z}>_{\mathbb{L}^{2}(m)}}F(Z_{u_{j}}\overline{Z}_{u_{j}})=E\otimes\mathbb{\mu
}_{x,y}(e^{\frac{1}{2}<\frac{L-\widehat{L}}{2}Z,\overline{Z}>_{\mathbb{L}%
^{2}(m)}}F(l^{u_{j}}+Z_{u_{j}}\overline{Z}_{u_{j}}))
\]

and the formula finally extends to

\begin{proposition}
For any bounded measurable function of a real field on $X:$%
\[
(\ast bis)\;E(Z_{x}\overline{Z}_{y}e^{\frac{1}{2}<\frac{L-\widehat{L}}%
{2}Z,\overline{Z}>_{\mathbb{L}^{2}(m)}}F(Z\overline{Z}))=E\otimes\mathbb{\mu
}_{x,y}(e^{\frac{1}{2}<\frac{L-\widehat{L}}{2}Z,\overline{Z}>_{\mathbb{L}%
^{2}(m)}}F(l+Z\overline{Z}))
\]
\end{proposition}

It induces the formula (*) obtained in the finite case if we consider the
trace of the process on any finite subset.

We see also that the restriction of the twisted Gaussian measure
\[
e^{\frac{1}{2}<\frac{L-\widehat{L}}{2}Z,\overline{Z}>_{_{\mathbb{L}^{2}(m)}}%
}P(dZ)
\]
to $\sigma(Z_{u}\overline{Z}_{u},u\in X)$ is a probability measure $Q$\ under
which the distribution of the process $(Z_{u}\overline{Z}_{u},u\in X)$ is
infinitely divisible. Moreover, the important fact is that this probability is
absolutely continuous with respect to the restriction of $P$ to $\sigma
(Z_{u}\overline{Z}_{u},u\in X)$. It is clear that that formula (**) of the
corollary extends in the same way to yield the following

\begin{corollary}
For any bounded measurable function of a nonnegative field on $X:$%
\[
(\ast\ast bis)\;\int\rho_{x}F(\rho_{u},u\in X)Q(d\rho)=\int\mathbb{\mu}%
_{x,x}(F(l^{u}+\rho_{u},u\in X))Q(d\rho)
\]
\end{corollary}

Hence it follows for example that the continuity of the Gaussian field $Z$
implies the continuity of the local time field $l$ under all loop measures
$\mu_{xx}$.

\bigskip

Note finally that these results, as in the symmetric case, can be extended to
some situations where the local time does not exist (like the two dimensional
Brownian motion), by considering the centered occupation field\ and the ''Wick
square'' $:Z_{x}\overline{Z}_{x}:$ (formally given by $Z_{x}\overline{Z}%
_{x}-K(x,x)$)\ as generalized random fields. This makes sense for the Wick
square provided $K$ is a Hilbert-Schmidt operator.

\end{document}